\documentclass[]{amsart}
\usepackage{amsmath,amssymb,color,hyperref}
\def \B{\mathcal{B}}
\def\marg#1{}

\theoremstyle{definition}

\newtheorem*{metaquestion}{Meta-question}
\newtheorem*{detailedquestions}{Detailed questions}

\title[Garside groups and geometry -- in memory of Patrick Dehornoy]{Garside groups and geometry}
\author{Bert Wiest}
\date{}

\begin{document}

\maketitle\marg{The title for the published version is ``Garside groups and geometry: some thoughts for Patrick Dehornoy''. For arXiv, it is ``Garside groups and geometry''}

\begin{abstract}
This article in memory of Patrick Dehornoy (1952 -- 2019) is an invitation to Garside theory for mainstream geometric group theorists interested in mapping class groups, curve complexes, and the geometry of Artin--Tits groups.
\end{abstract}


\section{Meeting Patrick}\label{S:Introduction}

When I remember Patrick Dehornoy, the first thing that comes to my mind is his humour -- he was one of the 
people who knew how to do serious work and have serious discussions with a smile and even a laugh, always in a way that advanced the work.  And his work! His energy and his efficiency at everything he touched were downright scary to me. He was clearly the driving force behind the two book projects \cite{DDRW1, DDRW2} he had with me (together with Ivan Dynnikov and Dale Rolfsen). Each time, his plans for these books grew more and more ambitious even during the writing process, and it must have been terribly frustrating to him that his co-authors had trouble following his rhythm. His unbelievable productiveness and efficiency concerned not only his mathematical work, but also his administrative work (as head of his department, but also in French national scientific administration), and his private projects (constructing his house and his garden, and creating humorous little movies, for instance). 

I remember Patrick's sometimes startling openness when discussing his private life. 
I remember his kindness and generosity, in particular with younger mathematicians, taking time and sharing his ideas. 
I remember his patience with less efficient people.
When I told him quite frankly that the book on the foundations of Garside theory was too long and abstract for me, he listened patiently and promised to write something short and enjoyable that would at least give a glimpse of the book. The result was the survey paper \cite{DehornoyQuadrNorm}, which I highly recommend as readable introduction to the subject. 

He must have known\marg{very few people $\to$ a very few experts} 
that some of his work was only readable for a very few experts. But I believe he had a clear vision of the kind of larger theory his work will one day be part of, and he had the confidence, and self-confidence, that it will be recognised as important. 
There is a scene in one of Patrick's\marg{added reference} little movies~\cite{DehornoyParadis} where a mathematician in the far future (played by C\'edric Villani) bases his revolutionary work on Patrick's LD-systems.  This was of course a joke, but who believes that it was \emph{only} a joke? Here is another telling example of his humour:
I once wrote to Patrick that a research project which we wanted to start had to wait, because I was still too busy with other projects. His answer: ``Quel dommage pour toi... dans ce cas tu n'auras pas la médaille Fields qui, pourtant, était à ta portée \ \ ;--)''.  (What a pity for you... you won't have the fields medal, which had been within your reach.)

\section{Garside theory}

Some of Patrick Dehornoy's major contributions in the last few years of his life were, of course, to the foundations of Garside theory~\cite{DDGKM}. 
The beautiful subject of Garside theory is now an important tool for geometric group theorists working on Artin-Tits groups -- see e.g.\  \cite{Bestvina, CMW, CharneyArtinBiautom,McCamEuclArtin, McCammSulway, BessisDual, BradyWattNoncross, AGJM2, MultifractRed1, MultifractRed2, MultifractRed3, MultifractRed4}. (And by the way, Patrick advocated the name ``Artin--Tits groups'' for this family of groups, in place of ``Artin groups''.)
On the other hand, Garside theory is not so much used by geometric group theorists working on mapping class groups, curve complexes, Teichmüller spaces etc. 
The rest of the present paper is meant as an invitation to them, and many others, to have a closer look at Garside theory.

This paper is mostly a list of open questions concerning Garside theory, questions whose answers should be of great interest to geometers. Even the most basic examples of Garside groups, namely the braid groups equipped with the classical Garside structure, are fascinating! All these questions are related to the following

\begin{metaquestion} Is Garside theory on braid groups a geometric theory?
\end{metaquestion}

We give substance to this vague question with three

\begin{detailedquestions} \textcolor{white}{ }
\begin{enumerate}
\item[(A)] Do Garside normal forms of braids represent quasi-geodesics in the curve graph?
\item[(B)] Is the additional length graph of braid groups, equipped with their classical Garside structure, quasi-isometric to the curve graph?
\item[(C)] Can Garside theory furnish a polynomial time solution to the conjugacy problem in braid groups?
\end{enumerate}
\end{detailedquestions}

In order to make Question~(A) more precise, we start with the following special case of an observation of  Masur and Minsky~\cite{MM1}:
the curve graph of the $n$-times punctured disk $\mathcal C(D_n)$ is quasi-isometric to an ``electrified'' Cayley graph of the braid group. Specifically, starting with $Cay(\B_n,\mathcal S_{Gar})$, the Cayley graph of~$\B_n$, equipped with the Garside generators, 
we can ``squash down'', or ``cone off'' the normaliser of every parabolic subgroup of the form $\langle \sigma_i,\ldots, \sigma_j\rangle$ with $1\leqslant i\leqslant j<n$, and also the subgroup $\langle\Delta^2\rangle$ generated by the Dehn twist  along the boundary of the disk. We do this simply by adding every element of every such subgroup to the generating set, thus adding plenty of edges to the Cayley graph; we denote $Cay(\B_n,\mathcal S_{\mathrm{NP}})$ the Cayley graph with respect to this enlarged, infinite generating set. The observation now is that $Cay(\B_n,\mathcal S_{\mathrm{NP}})$ is quasi-isometric to the curve graph $\mathcal C(D_n)$! 
Thus the composition\marg{q.i. $\to$ q.isom}
 $Cay(\B_n,\mathcal S_{Gar})\hookrightarrow Cay(\B_n,\mathcal S_{\mathrm{NP}})\stackrel{\mathrm{q.isom}}{\longrightarrow} \mathcal C(D_n)$ is a Lipschitz map from the standard Cayley graph to the curve graph.
We are now ready to give a more precise version of Question~(A):\marg{expanded (clarified) this sentence}
the images (or ``shadows'') of Garside normal form paths in $Cay(\B_n,\mathcal S_{Gar})$ under the above Lipschitz map are paths in~ $\mathcal C(D_n)$ -- do these paths form a uniform family of re-parametrised quasi-geodesics?

Why should we care about Question~(A)? One of the major open problems in the theory of mapping class groups is whether all mapping class groups are bi-automatic. And one of the most successful techniques used in the study of mapping class groups and related objects has been Masur--Minsky's hierarchy machinery on the curve graph -- see e.g.\ \cite{MM2,BehrstockHagenSisto2}. One should expect this technique to be helpful at attacking the bi-automaticity question.

A good place to start any such attack would surely be to understand what the \emph{known} bi-automatic normal forms on mapping class groups look like from the point of view of the curve graph.
But there is currently only one known bi-automatic normal form 
on any mapping class group, namely the Garside normal form (including the dual structure of~\cite{BKL}), on the braid groups.
The fact that Question~(A) is open is therefore a major embarrassment, and it illustrates the difficulty of the bi-automaticity question.
It is certainly reasonable to believe in a positive answer to Question~(A), but~\cite{RafiVerberne} is a cautionary tale.

Concerning Question~(B): in \cite{CalvezWiest1, CalvezWiest2} it was shown how to associate to every Garside group~$G$ a $\delta$-hyperbolic graph~$\mathcal C_{AL}(G)$ with an isometric $G$-action. It was shown that for all spherical-type Artin-Tits groups, the graph\marg{was ``complex''} 
is of infinite diameter, and the action has the weak proper discontinuity (``WPD'') property~\cite{BestvinaFujiwara}. 
The construction of~$\mathcal C_{AL}(G)$ starts again with the Cayley graph of~$G$ with respect to the Garside generators, and proceeds again by a ``coning off'', or ``adding to the generating set'' procedure: the idea is to add all the obvious elements -- obvious from the point of view of Garside theory -- which should be squashed so as to force hyperbolicity. 
Firstly, all powers of the Garside element~$\Delta$ must be added to the generating set (in order to squash the center of~$G$). Moreover, 
\marg{Changed the definition here, which was false, and also used the word ``braid''. The whole explanation of $C_{AL}$ is rewritten.}
if two elements $x$ and~$y$ can ``telescope into each other'',  
i.e.\ if $x$, $y$ and $x\cdot y$ have infimum zero (meaning that they can be written as products Garside generators, not their inverses, but that no such writing can use the letter~$\Delta$), and if all three  elements $x$, $y$, and $x\cdot y$ have the same length, then these three must belong to to generating set. (The element~$y$ in the preceding discussion is ``absorbable'', in the notation of \cite{CalvezWiest1, CalvezWiest2, CalvezWiest3}.) The obvious example of this kind telescoping takes place in the braid groups, 
when $x$ and~$y$ are positive braids of the same length with support in disjoint subsurfaces.

The proof that $\mathcal C_{AL}(G)$ is hyperbolic is quite easy from the construction. 
In order to give some idea of the proof that it has infinite diameter, we can state a corollary of this proof (which is unfortunately not stated in \cite{CalvezWiest2}): there are elements in the group whose axes in the Cayley graph have the strong contraction property, or equivalently, the bounded geodesic image property -- see \cite[Section 2]{ArzhCashenTao} for a good introduction to this property. (One fact which is helpful in proving the bounded geodesic image property is that Garside normal form paths are \emph{geodesics} in the Cayley graphs of Garside groups, with respect to their Garside generators \cite{CharneyMeier}.)

Question~(B) amounts to the following:
if two braids $x$ and $y$ can ``telescope into each other''  as above,
is there necessarily something in the geometry, in the subsurface structure of these braids, which explains this behaviour, like there is in the ``obvious example''?

Why should we care about Question~(B)? 
A positive answer would deepen our understanding both of Garside theory and of the curve complex. For a start, it would imply a positive answer to question~(A). Also, if $x$ and $y$ are two positive braids (with infimum zero), 
there would be a wonderfully elegant answer to the question: where is the quasi-center of a geodesic triangle in the curve graph (the electrified braid group) with vertices $1$, $x$ and $y$? The answer would be: at the greatest common divisor $x\wedge y$. Another attractive aspect of this question is its relation to a very puzzling gap in our understanding of the curve graph, namely: if you have a Euclidean disk of sufficiently large radius quasi-isometrically embedded in the Cayley graph of a mapping class group,  does electrifying the Cayley graph (projecting to the curve complex) shrink it to universally bounded diameter? As shown in \cite{CalvezWiest3}, a positive answer would imply a positive answer to Question~(B)!

When we look at more general Artin-Tits groups, we also see quite momentous consequences of a generalised version of Question~(B): 
there is currently a lot of activity trying to define an analogue of the curve graph for general Artin-Tits groups (a long-term goal being a proof that they are hierarchically hyperbolic \cite{BehrstockHagenSisto2}). There is actually an obvious candidate graph, which might conceivably work for every Artin-Tits group~$A$, namely the graph of irreducible parabolic subgroups~$\mathcal C_{parab}(A)$ defined in~\cite{CGGMW}. Here is the definition: for any Artin-Tits group~$A$, vertices of the graph $\mathcal C_{parab}(A)$ correspond to irreducible parabolic subgroups of~$A$, and two vertices are connected by an edge if, either, one of the corresponding subgroups is contained in the other, or if the two corresponding subgroups 
have trivial intersection and commute.  
In the case of braid groups, this graph is \emph{isometric} to the curve graph. However, in general, the hyperbolicity of these candidate graphs is an open problem, except in the cases $A_n, B_n, \widetilde A_n$ and $\widetilde C_n$, where it can be deduced (in a highly non-trivial way~\cite{CalvezCisnerosdelaCruz}) from the hyperbolicity of the curve graph. See \cite{CGGMW} for details on the case of Artin-Tits groups of spherical type, \cite{CalvezCisnerosdelaCruz} for the spherical and Euclidean type, and \cite{MorrisWright} for FC type.  

Getting back to the motivation for question~(B), the most elegant strategy for proving hyperbolicity of the graphs of irreducible parabolic subgroups, at least in the spherical-type case, and maybe also in the Euclidean-type case \cite{McCamEuclArtin,McCammSulway}, would be to prove that they are quasi-isometric to the additional length graphs of the corresponding groups, equipped with their Garside structures!

Turning to Question~(C), we recall that solving the conjugacy problem in mapping class groups in polynomial time is one of the major open problems concerning these groups, so any progress on Question~(C) would be remarkable.
A more precise version of Question~(C) is:  
``For any fixed number of strands~$n$, is there an algorithm which takes as its input two elements of the braid group~$\B_n$, of lengths $\ell_1$ and $\ell_2$, and which outputs the information whether or not the braids are conjugate, and which runs in time $P(\ell_1+\ell_2)$, for some polynomial~$P$?'' (Even more optimistically, we could ask for a polynomial dependance on both $\ell_1+\ell_2$ and the number of strands~$n$.) 

The biggest roadblock to solving Question~(C) using Garside theory is the following problem: consider a braid $x\in \B_n$ of length~$\ell$ which is \emph{rigid} (i.e.\ whose Garside normal form is \emph{cyclically} in normal form). We try to find an upper bound for the number of rigid conjugates of~$x$ -- the desired answer would be a polynomial $P_n(\ell)$, and ideally, the degree of $P_n$ should not grow too fast with respect to the number of strands~$n$. This is a purely Garside-theoretical question! Also, we know from Garside-theoretical considerations \cite{CarusoWiest} that \emph{generically}, $x$~has only $2\ell$ rigid conjugates, but we are interested in a worst-case bound. 
Now, computational experiments with large numbers of randomly generated braids, as well as explicit \marg{theoretical $\to$ explicit} constructions of families of braids with many rigid conjugates were reported on in~\cite{SchleimerWiest}. This work suggest that the desired upper bound may be of the form $|\{$rigid conjugates$\}|\leqslant c_n\cdot \ell^{n-2}$. Not only that, but it appears that when the set of rigid conjugates gets so exceptionally large, there is a geometric reason for this: the braid is ``almost reducible in multiple ways''. 
If a braid is reducible, then its closure admits a torus 
containing an ``interior braid'' which can slide freely relative to the ``exterior braid''. In our ``almost reducible braids'', by contrast, we do not have a torus but only a long tube whose two ends do not quite close up (what is called an ``ouroboros''). As a consequence, the interior of the tube can only slide by a limited amount relative to the exterior. If there are several such tubes in the same closed braid, which can all slide by some amount relative to each other, then the set of conjugates thus obtained can grow polynomially with the length~$\ell$.\marg{$l \to \ell$}

In particular, it appears that the braids for which solving the conjugacy problem using Garside theory is hard all have very short translation length in the curve complex. In terms of hierarchies, there are multiple disjoint or nested subsurfaces with large subsurface projections.

In summary, in each of the cases referenced by Questions (A), (B), and~(C), there may be some hidden rigidity to Garside theoretical structures, and in each case, this rigidity may correspond to some geometric information. If true, these phenomena may prove mutually beneficial for both geometry and Garside theory. And if not, well, then Garside may reveal structures not visible in curve-complex like objects.\marg{Added this last sentence.}
I imagine that either answer would have been deeply satisfying for Patrick Dehornoy.


\begin{thebibliography}{00}

\bibitem{AGJM2} {\bf Samy Abbes, Sébastien Gouëzel, Vincent Jugé, Jean Mairesse},  {\it Asymptotic combinatorics of Artin-Tits monoids and of some other monoids}, J.\ Algebra 525 (2019), 497--561

\bibitem{ArzhCashenTao} {\bf Gulnara N. Arzhantseva, Chris H. Cashen, Jing Tao}, {\it Growth tight actions}, Pacific J.\ Math.~278 (2015), no. 1, 1--49.

\bibitem{BehrstockHagenSisto2} {\bf Jason Behrstock, Mark F.\ Hagen, Alessandro Sisto} {\it Hierarchically hyperbolic spaces~II: Combination theorems and the distance formula}, arXiv:1509.00632

\bibitem{BessisDual} {\bf David Bessis}, {\it The dual braid monoid}, Ann.\ Sci.\ École Norm.\ Sup.(4), (2003), no. 36(5), 647--683

\bibitem{Bestvina} {\bf Mladen Bestvina}, {\it Non-positively curved aspects  of Artin groups of finite type}, Geometry and Topology 3 (1999),  269--302.

\bibitem{BestvinaFujiwara} {\bf Mladen Bestvina, Koji Fujiwara}, {\it Bounded cohomology of subgroups of mapping class groups}, Geom. Topol. 6 (2002), 69--89

\bibitem{BKL} {\bf Joan Birman, Ki Hyoung  Ko, Sang Jin  Lee}, {\it A new approach to the word and conjugacy problems in the braid groups}, Adv.\ Math.\ 139 (2), (1998) 322--353.

\bibitem{BradyWattNoncross} {\bf Thomas Brady, Colum Watt}, {\it Non-crossing partition lattices in finite real reflection groups}, Trans.\ Amer.\ Math.\ Soc.\ 360 (2008), no. 4, 1983--2005. 

\bibitem{CalvezCisnerosdelaCruz} {\bf Matthieu~Calvez, Bruno A. Cisneros de la Cruz}, {\it Curve graphs for Artin-Tits groups of type $B$, $\tilde A$ and $\tilde C$ are hyperbolic}, arXiv:2003.04796 

\bibitem{CalvezWiest1} {\bf Matthieu Calvez, Bert Wiest}, {\it Curve graphs and Garside groups}, Geometriae Dedicata 188(1) (2017), 195--213

\bibitem{CalvezWiest2} {\bf Matthieu Calvez, Bert Wiest}, {\it Acylindrical hyperbolicity and Artin-Tits groups of spherical type}, Geometriae Dedicata 191(1) (2017), 199--215

\bibitem{CalvezWiest3} {\bf Matthieu Calvez, Bert Wiest}, {\it Hyperbolic structures for Artin-Tits groups of spherical type},
arXiv:1904.02234, to appear in Advances in Mathematics

\bibitem{CarusoWiest} {\bf Sandrine Caruso, Bert Wiest}, {\it On the genericity of pseudo-Anosov braids II: conjugations to rigid braids}, Groups Geom.\ Dyn.~11 (2017), 549--565.

\bibitem{CharneyArtinBiautom} {\bf Ruth Charney}, {\it Artin groups of finite type are biautomatic}, Math. Ann. 292 (1992), no. 4, 671--683. 

\bibitem{CharneyMeier} {\bf Ruth Charney, John Meier}, {\it The language of geodesics for Garside groups}, Math.\ Zeitschrift, 248 (2004), 495--509.

\bibitem{CMW} {\bf R.~Charney, J.~Meier, K.~Whittlesey}, {\it Bestvina's normal form complex and the homology of Garside groups},  Geom.\ Dedicata 105 (2004), 171--188.

\bibitem{CGGMW} {\bf Mar\'ia Cumplido, Volker Gebhardt, Juan Gonz\'alez-Meneses, Bert Wiest}, {\it On parabolic subgroups of Artin-Tits groups of spherical type}, Advances in Mathematics 352 (2019), 572--610

\bibitem{DehornoyParadis}\marg{Added this reference} {\bf Patrick Dehornoy}, {\it Le Paradis des mathématiciens}, short film, https:// vimeo.com/205778279

\bibitem{DehornoyQuadrNorm} {\bf Patrick Dehornoy}, {\it Garside and quadratic normalisation: a survey}, arXiv 1504.07788, Proceedings DLT 2015, I.Potapov ed., Springer LNCS 9138, pp. 14-45, 

\bibitem{MultifractRed1}{\bf Patrick Dehornoy}, {\it Multifraction reduction I: The 3-Ore case and Artin-Tits groups of type~FC}, J.~Comb.\ Algebra~1 (2017) 185--228

\bibitem{MultifractRed2}{\bf Patrick Dehornoy}, {\it Multifraction reduction II: Conjectures for Artin-Tits groups}, J.~Comb.\ Algebra~1 (2017) 229--287

\bibitem{DDGKM} {\bf Patrick Dehornoy, Fran\c{c}ois Digne, Eddy Godelle, Daan Krammer, Jean Michel}, {\it Foundations of Garside Theory}, EMS Tracts in Mathematics, volume 22, European Mathematical Society, 2015.

\bibitem{DDRW1} {\bf Patrick Dehornoy, with Ivan Dynnikov, Dale Rolfsen, Bert Wiest}, {\it Why are braids orderable?}, Panoramas et synthèses n.~14, Société Mathématique de France (2002)

\bibitem{DDRW2} {\bf Patrick Dehornoy, with Ivan Dynnikov, Dale Rolfsen, Bert Wiest}, {\it Ordering braids}, Surveys and Monographs vol.~148, American Mathematical Society (2008)

\bibitem{MultifractRed4}{\bf Patrick Dehornoy, Derek Holt, Sarah Rees}, {\it Multifraction reduction IV: Padding and Artin-Tits groups of sufficiently large type}, J.~Pure Appl.\ Algebra 222 (2018) 4082--4098

\bibitem{MultifractRed3}{\bf Patrick Dehornoy, Friedrich Wehrung}, {\it Multifraction reduction III: The case of interval monoids}, J.~Comb. Algebra~1 (2017) 341--370

\bibitem{MM1} {\bf Howard Masur, Yair Minsky}, {\it Geometry of the complex of curves I: Hyperbolicity}, Invent. math. 138 (1999), 103--149.

\bibitem{MM2} {\bf Howard Masur, Yair Minsky}, {\it Geometry of the complex of curves II: Hierarchical Structure}, Geometric and Functional Analysis, 10(4)(2000), 902--974 .

\bibitem{McCamEuclArtin} {\bf Jon McCammond}, {\it The structure of Euclidean Artin groups},  Geometric and cohomological group theory (London Math.\ Soc.\ Lecture Note Series), 
Cambridge University Press, (2017), 82--114.

\bibitem{McCammSulway} {\bf Jon McCammond, Richard Sulway}, {\it Artin groups of Euclidean type}, Invent. Math. 210, 2017, 231--282.

\bibitem{MorrisWright} {\bf Rose Morris-Wright}, {\it Parabolic subgroups of Artin groups of FC type}, arXiv:1906.07058

\bibitem{OsinAcylHyp} {\bf Denis Osin}, \textit{Acylindrically hyperbolic groups}, Trans. Amer. Math. Soc. 368 (2016), 851-888.

\bibitem{RafiVerberne} {\bf Kasra Rafi, Yvon Verberne}, {\it Geodesics in the mapping class group}, arXiv:1810.12489

\bibitem{SchleimerWiest} {\bf Saul Schleimer, Bert Wiest}, {\it Garside theory and subsurfaces: some examples in braid groups}, Groups Complexity Cryptology, 11 (2019), no. 2, 61–75.

\end{thebibliography}
\end{document}